\magnification 1200
\parindent=0mm
\let\dis=\displaystyle
\let\ge=\varepsilon
\let\ra=\rightarrow  
\def\F{{\cal F}}
\def\R{{\bf R}} \def\E{I\!\!E}

\long\def\titre#1{\vbox{\leftskip=0pt plus 1fil\rightskip=0pt plus 1fil
\parfillskip=0pt#1}}
\long\def\abstract#1{\centerline{\hsize12cm \vbox{{\bf Abstract\kern 2pt.}{ #1}}}}
\long\def\resume#1{\centerline{\hsize12cm \vbox{{\bf R\'esum\'e\kern 2pt.}{ #1}}}}
\def\fin{\kern 6pt\hbox{\vrule height 4pt  width 1mm}}

\titre{{\bf Stochastic differential equations with non-Lipschitz coefficients:\par
II. Dependence with respect to initial values}
\medskip
Shizan FANG and Tusheng ZHANG}
\vskip 10mm
\medskip
{\sevenrm S.F.: I.M.B, UFR Sciences et techniques, Universit\'e de Bourgogne, 9 avenue Alain Savary, B.P. 47870, 21078 Dijon, France.}
\smallskip
{\sevenrm T.Z.: Department of Mathematics, University of Manchester, Oxford road, Manchester, M13 9PL, England.}
\vskip 10mm

\abstract{The existence of the unique strong solution for a class of stochastic differential equations with non-Lipschitz coefficients was established recently. In this paper, we shall investigate the dependence with respect to the initial values. We shall prove that the non confluence of solutions holds under our general conditions. To obtain a continuous version, the modulus of continuity of coefficients is assumed to  be less than 
$\dis |x-y|\log{1\over|x-y|}$. In this case, it will give rise to a flow of homeomorphisms if the coefficients are compactly supported.}
\vskip 10mm

{\bf 0. Introduction}
\medskip
\quad Let $ \sigma:\ R^d\rightarrow  R^d\otimes R^m$ and 
$b:\ R^d\rightarrow R^d$ be continuous functions. It is well known (see [IW], [RY]) that the following It\^o s.d.e:
$$ dX(t)=\sigma(X(t))\,dW_t + b(X(t))\, dt,\quad X(0)=x_o\leqno(0.1)$$

has a weak solution up to a lifetime $\zeta$. It has been proved recently in [FZ2] (see [FZ1] for a short version) that the s.d.e. $(0.1)$ has a unique strong solution if the coefficients $\sigma$ and $b$ satisfy the following non-Lipschitz conditions

$$\dis
\left\{\matrix{||\sigma(x)-\sigma(y)||^2&\leq& C\, |x-y|^2\,r(|x-y|^2),\cr
|b(x)-b(y)|&\leq& C\,  |x-y|\,r(|x-y|^2)}\right.
\leqno(H1)$$

for $|x-y|\leq \delta_o$, where $r: ]0,\delta_o]\ra \R_+$ is a strictly positive function of class $C^1$ satisfying
 (i)\ $\dis\lim_{s\ra 0} r(s)=+\infty$,
\ (ii)\ $\dis\int_0^{\delta_o} {ds\over sr(s)}=+\infty$  and
\ (iii)\ $\dis\lim_{s\ra 0}{sr'(s)\over r(s)}=0$.
\medskip

\quad Let $\dis X_t(x_o,w)$ be the solution of the s.d.e. $(0.1)$. It is well-known that the solution admits  a continuous version $\dis\tilde X_t(x_o,w)$ if the coefficients $\sigma$ and $b$ are locally Lipschitzian (see [Pr]), and it gives rise to a flow of homeomorphisms if the coefficients are globally Lipschitzian (see [Ku]). We refer also to [Ma2], [IW] and [El] for the study of stochastic flow of diffeomorphisms, to [Ma1] for the non-Lipschitz feature in the study of homeomorphisms of the circle $S^1$. Many interesting phenomena for stochastic differential equations with  non-Lipschitz coefficients   have been elucidated in [LJR1,2]. In this work, we shall investigate the dependence of solutions with respect to the initial values under the non-Lipschitz condition $(H1)$, which can not be covered in  [LJR1,2] and generalize [Ma1] to general situation.
\medskip

\quad The organization of the paper is as follows. In section 1, we shall prove that the non confluence (or non-contact) property holds under $(H1)$. The function $r$ considered in $(H1)$ includes obviously all functions like $\dis\log{1\over\xi},\, \log{1\over\xi}\,\log{\log{1\over\xi}}, \cdots$. Such kind of properties were also studied by 
M. Emery in an early work [Em] for Lipschitz case, and by T.Yamada and Y. Ogura for a non-Lipschitz case in [YO]. The conditions in [YO] includes  also the function $\dis\xi\log{1\over\xi}$, but their mixed condition about $\sigma$ and $b$
$$\dis \lim_{s\downarrow 0}\,{\sup_{s\leq t\leq \delta_o}\Bigl[\kappa(t)\int_t^{\delta_o}{du/\rho^2(u)}\Bigr]
\over \int_s^{\delta_o}\!\!\Bigl[\int_t^{\delta_o}{du/ \rho^2(u)}\Bigr]dt}=0$$
where $\dis \rho^2(u)=u^2r(u^2)$ and $\kappa(u)=ur(u^2)$, is not easy  to be checked in general. In section 2, we shall establish a continuous version of the solutions. The main tool for doing this is the Kolmogorov's  modification theorem, which will work in the case where $\dis r(\xi)=\log{1\over\xi}$. However, it seems difficult to apply  the modification theorem in  the case where 
$\dis r(\xi)=\log{1\over\xi}\,\log{\log{1\over\xi}}$. Finally in section 3, we shall prove that  the continuous version of solutions will give rise to a flow of homeomorphisms if the coefficients are compactly supported.  
\vskip 8mm

{\bf 1. Non contact property}
\bigskip

{\bf Theorem 1.1}\quad {\it Assume (H1) and the s.d.e $(0.1)$ has no explosion. If  $x_o\neq y_o$, then  $X_t(x_o,w)\neq X_t(y_o,w)$  almost surely  for all $t>0$.}

\medskip
{\bf Proof.}
Consider
$$ \psi(\xi) = \int_{\xi}^1 {ds\over sr(s)}
\quad\hbox{\rm and}\quad \Phi(\xi)= e^{\psi(\xi)},\quad\hbox{\rm for}\ 1\geq \xi> 0. $$
We have

$$\dis \Phi'(\xi)=-{\Phi(\xi)\over \,\xi\,r(\xi)}\leq 0.$$
By condition (i) and (iii) about the function $r$, it exists $1>\delta>0$ and a constant $C>0$ such that

$$\dis \Phi''(\xi) =
{\Phi(\xi)\,\bigl(1+r(\xi)+{\xi}r'({\xi})\bigr)\over (\xi r({\xi}))^2}\leq C \Phi(\xi){1\over \xi^2r(\xi)}
\leqno(1.1)$$
for $0<\xi<\delta$. To be simplified, denote $X_t(x_o)=X_t(x_o,w)$. 
Let  $\dis\eta_t= X_t(x_o)-X_t(y_o)$ and
$\dis \xi_t(w)=|\eta_t(w)|^2$. We have

$$\dis\eqalign{
d\xi_t=&2\bigl<\eta_t, \bigl(\sigma(X_t(x_o))-\sigma(X_t(y_o))\bigr)\,dW_t\bigr>\cr
&+ 2\bigl<\eta_t, b(X_t(x_o))-b(X_t(y_o))\bigr>\,dt\cr
& +||\sigma(X_t(x_o))-\sigma(X_t(y_o))||^2\, dt,}
$$

and the stochastic contraction $d\xi_t\cdot d\xi_t$ is given by
$$ 4|\bigl(\sigma^*(X_t(x_o))-\sigma^*(X_t(y_o))\bigr)\eta_t|^2\,dt$$
where $\sigma^*$ denotes the transpose matrix of $\sigma$.
Without loss of generality, we may assume $|x_o-y_o|<{\delta\over 2}$. 
Let $\dis 0<\ge<|x_o-y_o|<{\delta\over 2}$. Define
$$\dis \tau_\ge=\inf\bigl\{t>0, \xi_t\leq\ge\bigr\},
\quad \tau=\inf\bigl\{t>0, \xi_t=0\bigr\}.$$

We have $\dis \tau_\ge\uparrow \tau$ as $\ge\downarrow 0$.
Let $$ \dis \zeta =\inf{\bigl\{t>0, \ \xi_t\geq {3\delta\over 4}\,\bigr\}}.$$

Using It\^o formula, 
$$\dis\eqalign{\Phi\bigl(\xi_{t\wedge\tau_\ge\wedge\zeta}\bigr)
&=\Phi(\xi_o) + 2\int_0^{t\wedge\tau_\ge\wedge\zeta}
\Phi'({ \xi_s})
\bigl<\eta_s, \bigl(\sigma(X_s(x_o))-\sigma(X_s(y_o))\bigr)\,dW_s\bigr>\cr
&\hskip 3mm
+2\int_0^{t\wedge\tau_\ge\wedge\zeta}
\Phi'({ \xi_s})
\bigl<\eta_s, b(X_s(x_o))-b(X_s(y_o))\bigr>\,ds\cr
&\hskip 3mm
+\int_0^{t\wedge\tau_\ge\wedge\zeta}
\Phi'({ \xi_s})
||\sigma(X_s(x_o))-\sigma(X_s(y_o))||^2\, ds\cr
&\hskip 3mm
+2\int_0^{t\wedge\tau_\ge\wedge\zeta}
\Phi''({\xi_s})
|\bigl(\sigma^*(X_s(x_o))-\sigma^*(X_s(y_o))\bigr)\eta_s|^2\,ds.}\leqno (1.2)$$

By $(H1)$ we have
$$|\bigl(\sigma^*(X_s(x_o))-\sigma^*(X_s(y_o))\bigr)\eta_s|^2\leq \xi_s^2r(\xi_s).$$ 
Combining with (1.1), we get
$$\dis 2\Phi''({\xi_s})
|\bigl(\sigma^*(X_s(x_o))-\sigma^*(X_s(y_o))\bigr)\eta_s|^2
\leq C\, \Phi(\xi_s).$$

Again by $(H1)$ and the expression of $\Phi'(\xi)$, we have
$$\dis \Bigl|\Phi'({ \xi_s})
\bigl<\eta_s, b(X_s(x_o))-b(X_s(y_o))\bigr>\Bigr|
\leq C\, \Phi(\xi_s).$$

Therefore according to $(1.2)$,
$$\dis\Phi\bigl(\xi_{t\wedge\tau_\ge\wedge\zeta}\bigr)
\leq\Phi(\xi_o) + \hbox{\rm martingale}
+ C\, \int_0^{t\wedge\tau_\ge\wedge\zeta} \Phi(\xi_s)\, ds.$$

Taking expectation, we get

$$\dis\eqalign{ \E\Bigl(\Phi(\xi_{t\wedge\tau_\ge\wedge\zeta})\Bigr)
&\leq \Phi(\xi_o)+ C\E\Bigl(\int_0^{t\wedge\tau_\ge\wedge\zeta} \Phi({\xi_s})ds\Bigr)\cr
&\leq \Phi(\xi_o)
+ C\,\int_0^t \E\bigl(\Phi(\xi_{s\wedge\tau_\ge\wedge\zeta})\bigr)\, ds}$$
which implies that

$$\dis\E\Bigl(\Phi(\xi_{t\wedge\tau_\ge\wedge\zeta})\Bigr)
\leq \Phi(\xi_o)e^{C\, t},\quad\hbox{\rm for all  }\ t>0.$$
Consequently,
$$P(\tau_\ge<t\wedge\zeta )\Phi(\ge )\leq \Phi(\xi_o)e^{C\, t},\quad\hbox{\rm for all  }\ t>0.$$
Now first letting $\ge \rightarrow 0$, we obtain 
$$P(\tau <t\wedge\zeta )=0,\quad\hbox{\rm for all  }\ t>0,$$
and then letting  $t\rightarrow \infty$ we get $P(\tau <\zeta )=0$. Therefore, $\xi_.$ is positive almost surely on the interval $[0,\zeta ]$.
Now define $T_0:=0$,
$$\dis T_1:=\zeta,\quad \quad T_2=\inf\bigl\{t>0,\ \xi_t\leq {\delta^2\over 4}\bigr\}$$
and generally 
$$\dis T_{2n}=\inf\bigl\{t>T_{2n-1},\ \xi_t\leq {\delta^2\over 4}\bigr\},\quad \quad T_{2n+1}=\inf\bigl\{t>T_{2n},\ \xi_t\geq {3\delta \over 4}\bigr\}$$
Clearly $T_n\rightarrow \infty$ almost surely as $n\rightarrow \infty$. By definition, $\xi_.$ is positive on the interval $[T_{2n-1}, T_{2n}]$. By pathwise uniqueness of solutions, $X$ enjoys the strong Markovian property. Starting again   from $T_{2n}$ and applying the same arguments as in the first part of the proof, one  can  show that $\xi_.$ is positive almost surely also on the interval $[T_{2n}, T_{2n+1}]$. This completes the proof.\fin

\vskip 0.3cm
{\bf Theorem 1.2}\quad {\it Assume that for $|x|\geq 1$,
$$\dis
\left\{\matrix{||\sigma(x)||^2&\leq& C\,( |x|^2\,\rho(|x|^2)+1),\cr
|b(x)|&\leq& C\,  (|x|\,\rho(|x|^2)+1)}\right.
\leqno(H2)$$
where $\rho: [1, +\infty[\ra \R_+$ is a function of class ${\cal C}^1$ satisfying
 $(i)\ \dis\int_1^\infty {ds\over s\rho(s)+1}=+\infty$,\ 
$(ii) \ \dis\lim_{s\ra+\infty}{s\rho'(s)\over \rho(s)}=0$ and
 $\ (iii)\ \dis\lim_{s\ra +\infty}\rho(s)=+\infty$.\ 
Then $\ \lim_{|x_o|\ra +\infty}|X_t(x_o,w)|=+\infty$ in probability.}

\medskip
{\bf Proof.}   Let $f\in C^1(R_{+})$ be a fixed,  strictly positive  $C^1$ function on $R_{+}$ that satisfies
$$\dis f(s)=\rho (s) \quad \quad \hbox{\rm for } \quad s\geq 1.$$
Then it is easy to see that $(H2)$ holds for all $x\in R^d$ , with $\rho$ replaced by $f$.
From now on, we will use $C$ to denote a generic constant which may change from line to line.

Define
$$ \dis \psi(\xi) = \int_0^\xi {{ds}\over {sf(s)+1}},\quad \xi\geq 0$$
and put
$$\dis  \Phi(\xi)=e^{-\psi(\xi) }.$$ 
Keeping the assumptions on $\rho$ in mind it follows that
$$\dis
 \Phi^{\prime}(\xi)=- {{\Phi(\xi)}\over {\xi f(\xi)+1}}\leq 0,$$
and
$$\dis\eqalign{
\Phi^{\prime\prime }(\xi)&= {{\Phi(\xi)}\over {(\xi f(\xi)+1)^2}}[1+f(\xi )+\xi f^{\prime}(\xi)]
\cr
&\leq C\Phi(\xi){{f(\xi )}\over {(\xi f(\xi )+1)^2}}.}
\leqno (1.3)$$

Let $ \xi (t)=|X_t(x_0)|^2$. For any constant $R>0$, define
$$\dis \tau_R =\inf{\big\{t\geq 0, \ |X_t(x_0)|\leq  R\,\big\}}.$$
By It\^o formula, we have
$$\dis\eqalign{
\Phi(\xi_{t\wedge\tau_R})
= &\Phi (|x_0|^2) + 2 \int_0^{t\wedge\tau_R} \Phi^{\prime}(\xi_s)
\big <X_s, \sigma(X_s)dW_s\big >\cr
+& 2\int_0^{t\wedge\tau_R} \Phi^{\prime}(\xi_s)\,
\big <X_s, b(X_s)\big >\,ds
+\int_0^{t\wedge\tau_R}\Phi^{\prime}(\xi_s)
||\sigma(X_s)||^2\,ds\cr
+&2  \int_0^{t\wedge\tau_R} \Phi^{\prime\prime }(\xi_s)\,|\sigma^*(X_s)X_s|^2\,ds.}
\leqno (1.4)
$$
By $(H2)$, 
it holds that
$$\dis  {|\sigma^*(X_s)X_s|^2\over (\xi_s\, f(\xi_s)+1)^2 }
\leq C\ {\xi_s\, (\xi_s\, f(\xi_s)+1)\over (\xi_s\, f(\xi_s)+1)^2 }$$
Together with $(1.3)$, we have
$$\dis
 \int_0^{t\wedge\tau_R}\Phi^{\prime\prime }(\xi_s)\,|\sigma^*(X_s)X_s|^2\,ds
\leq C \int_0^{t\wedge\tau_R} \Phi(\xi_s)\, ds.
$$
Similarly, we have for some constant $C>0$,
$$\dis
 \bigl|\Phi'(\xi_s)\big <X_s, b(X_s)\big >\bigr| \leq C\,\Phi(\xi_s), \quad s>0.
$$
Combining above inequalities, we get from (1.4)
$$\dis E\big (\Phi(\xi_{t\wedge\tau_R})\big )
\leq \Phi (|x_0|^2)+ C \int_0^t E\big (\Phi(\xi_{s\wedge\tau_R})\big )\, ds,$$
which implies that
$$\dis E\big (\Phi(\xi_{t\wedge\tau_R})\big )
\leq \Phi (|x_0|^2)\, e^{C t}.$$
This gives  that
$$ \dis P\big(\tau_R\leq t \big) \Phi (R^2)\leq E\big (\Phi(\xi_{t\wedge\tau_R})\big )
\leq \dis \Phi (|x_0|^2)\, e^{C t}.$$
Therefore,
$$\dis
 P\big(\inf_{0\leq s\leq t}
	|X_s(x_o)| \leq  R \big)
\leq  e^{C t}\, \exp\Bigl\{-\int_{R^2}^{|x_o|^2}{ds\over sf(s)+1}\Bigr\},$$
which tends to $0$ when $|x_o|\ra +\infty$.\fin

\vskip 8mm

{\bf 2. Continuous dependence of initial values}
\bigskip

\quad In this section, we shall show that the solution to (0.1) admits a version that  is jointly  continuous in $(t,x_o)$. Let's begin with the following lemma.
\medskip

{\bf Lemma 2.1}\quad {\it Let $p\geq 1$. Assume that the coefficients $\sigma$ and $b$ are compactly supported, say,
$$\dis \sigma(x)=0\quad\hbox{\rm and}\quad b(x)=0\quad\hbox{\rm for}\ |x|\geq R,
\leqno(2.1)$$

and the function $r$ in (H1) $ : ]0,\delta_o]\ra \R_+$ is decreasing. Then there exists a constant $C_p>0$ such that for $|x|\leq R+1$ and $|y|\leq R+1$,

$$\dis
\left\{\matrix{||\sigma(x)-\sigma(y)||^2
&\leq& C_p\, |x-y|^2\,r\Bigl(({|x-y|^2\over M})^p\Bigr),\cr
|b(x)-b(y)|&\leq& C_p\,  |x-y|\,r\Bigl(({|x-y|^2\over M})^p\Bigr)}\right.
\leqno(2.2)$$
where $\dis M={4(R+1)^2\over\delta_o}\geq 1$.}

\medskip
{\bf Proof.}\quad Because of the similarity, we only prove the conclusion for $b$. 
Let $\dis\delta=\inf\{{\delta_o\over 2}, {1\over 2}\}$. If  $\dis |x-y|\leq {\delta}$, by hypothesis $(H1)$,
$$\dis|b(x)-b(y)|\leq C\ |x-y|r(|x-y|^2)
\leq |x-y|\,r\Bigl(({|x-y|^2\over M})^p\Bigr),\leqno (2.3)$$
since $r$ is supposed to be decreasing. Remark that
$$\dis \inf_{{\delta}\leq\xi\leq 2(R+1)}\xi\, r\Bigl(\Bigl[{\xi^2\over M}\Bigr]^p\Bigr)
\geq {\delta }\,r\bigl({\delta}^p\bigr)$$
and
$\dis \sup_{x,y}|b(x)-b(y)|\leq 2 ||b||_\infty$. Therefore there exists a constant $C_p>0$ such that

$$\dis |b(x)-b(y)|
\leq C_p\, |x-y|\,r\Bigl(({|x-y|^2\over M})^p\Bigr)
\quad\hbox{\rm for }\ |x-y|\geq\delta.\leqno (2.4)$$

Combining (2.3) and (2.4), we get the result.\fin
\medskip
\medskip
{\bf Lemma 2.2}\quad {\it Assume the same hypothesis as in lemma 2.1 and furthermore 
$ \xi\ra \xi r(\xi)$ is concave over $ ]0, \delta_o]$.  
Let $p\geq 1$ be an integer. For $|x_o|\leq R+1$ and $|y_o|\leq R+1$, set
$$\dis \eta_t =|X_t(x_o)-X_t(y_o)|^2\quad\hbox{\rm and}
\quad \xi_t=({\eta_t\over M})^p$$

where $M$ is the constant appeared in (2.2). Put  $\phi (t)=E(\xi_t)$. Then we have  
$$\dis \varphi'(t)\leq C_p\, \varphi(t)\, r(\varphi(t))
\leqno(2.5)$$ 
for some constant $C_p$.}
\medskip
{\bf Proof.}\quad We remark that under the assumptions $(2.1)$, for any $|x_o|\leq R+1$, $|X_t(x_o)|\leq R+1$ almost surely for all $t\geq 0$.
In fact, this can be seen as follows. Define

$$\dis T=\inf\bigl \{t\geq 0;\ |X_t(x_o)|\geq R+1\bigr\}$$
Set  $\dis Y_t(x_0)=X_{t\wedge T}(x_o) $. Then
$$\dis Y_t=x_0+\int_0^{t\wedge T} \sigma (Y_s) dW_s+\int_0^{t\wedge T}b(Y_s) ds.$$
Since 
$$\dis \E\Bigl(\int_0^t \sigma^2(Y_s)\bigl( {\bf 1}_{(s<T)}-1\bigr)^2\, ds\Bigr)
=\E\Bigl(\int_T^t \sigma^2(Y_s)\, ds\Bigr)=0,$$
we have for $t\geq 0$,
$$\dis \int_0^{t\wedge T} \sigma (Y_s) dW_s=\int_0^{t} \sigma (Y_s) dW_s \quad\hbox{\rm and}\ \int_0^{t\wedge T}b(Y_s) ds=\int_0^{t}b(Y_s) ds,$$
almost surely. 
We see that $\{Y_t,\ t\geq 0\}$ satisfies the same stochastic differential equation as 
$\{X_t,\ t\geq 0\}$. By the pathwise uniqueness in [FZ2], we conclude that $Y_t=X_t$ a.s. for all $t\geq 0$, which proves the claim.
\medskip 
Now we shall proceed as in [Fa]. By It\^o's formula,
$$\dis\eqalign{
d\eta_t=&2\bigl<X_t(x_o)-X_t(y_o), 
\bigl(\sigma(X_t(x_o))-\sigma(X_t(y_o))\bigr)\,dW_t\bigr>\cr
&+ 2\bigl<X_t(x_o)-X_t(y_o), b(X_t(x_o))-b(X_t(y_o))\bigr>\,dt\cr
& +||\sigma(X_t(x_o))-\sigma(X_t(y_o))||^2\, dt,}$$

and
$$\dis d\xi_t={1\over M^p}\Bigl( p\eta_t^{p-1}d\eta_t 
+{1\over 2}p(p-1)\eta_t^{p-2}d\eta_t\cdot d\eta_t\Bigr).\leqno (2.6)$$

\medskip

\quad Now using $(2.2)$,
$$\dis \eqalign{&{2p\over M^p}\eta_t^{p-1}\Bigl|\bigl<X_t(x_o)-X_t(y_o),
 b(X_t(x_o))-b(X_t(y_0))\bigr>\Bigr|\cr
&\leq C_p {2p\over M^p}\eta_t^{p-1}\eta_t \, r(\xi_t)
=2pC_p\xi_tr(\xi_t).}$$

Similarly , we get the same control for other terms in (2.6) except the martingale part. Now by expression of $d\xi_t$, we have
$$\dis \xi_{t+\ge}-\xi_t \leq M_{t+\ge}-M_t + C_p\int_t^{t+\ge} \xi_sr(\xi_s)\, ds$$
where $M_t$ is the martingale part of $\xi_t$. It follows that
$$\dis \E(\xi_{t+\ge}|\F_t)-\xi_t \leq C_p\E\Bigl( \int_t^{t+\ge}\xi_sr(\xi_s)ds|\F_t\Bigr)$$
where $\F_t$ is the natural filtration generated by $\{w(s);\ s\leq t\}$. 
Let
$\dis\varphi(t)=\E(\xi_t)$. Then
$$\dis \varphi(t+\ge)-\varphi(t)\leq C_p\, \int_t^{t+\ge}\E(\xi_s\, r(\xi_s))\,ds.$$
It follows that
$$\dis \varphi'(t)\leq C_p\E(\xi_t\, r(\xi_t)).$$

Since
$$\dis \xi\ra \xi r(\xi)\quad\hbox{\rm is concave over }\ ]0, {\delta}],$$
we have 
$$\dis \varphi'(t)\leq C_p\, \varphi(t)\, r(\varphi(t)).
\leqno(2.7)$$

\medskip

{\bf Theorem 2.3}\quad{\it  Assume (H1) and the s.d.e $(0.1)$ has no explosion. 
Consider  $\dis r(\xi)=\log{1\over\xi}$. Then there exists a version of 
$X_t(w,x_o)$ such that $ (t,x_o)\ra X_t(w,x_o)$  is continuous over $[0,+\infty[\times\R^d$
almost surely.}

\medskip
{\bf Proof.}\quad It has been proved in [FZ2] that the s.d.e $(0.1)$ has no explosion under the hypothesis $(H2)$. We split the proof into two steps.

{\bf Step 1}. Assume that  $\sigma$ and $b$ are compactly supported, say,
$$\dis \sigma(x)=0\quad\hbox{\rm and}\quad b(x)=0\quad\hbox{\rm for}\ |x|\geq R.$$
Let $\varphi$ be defined as in Lemma 2.2. Solving (2.7), 
we get $\dis\varphi(t)\leq (\varphi(0))^{e^{-C_p\,t}}$
or explicitly
$$\dis \E\Bigl( |X_t(x_o)-X_t(y_o)|^{2p}\Bigr)
\leq C_p\, |x_o-y_o|^{2pe^{-C_p\, t}}.$$

On the other hand, it is easy to see that
$$\dis \E\Bigl( |X_t(x_o)-X_s(x_o)|^{2p}\Bigr)
\leq C_p |t-s|^p.$$

Therefore,
$$\dis \E\Bigl( |X_t(x_o)-X_s(y_o)|^{2p}\Bigr)
\leq C_p\, \Bigl[|t-s|^p + |x_o-y_o|^{2pe^{-C_p\, t}}\Bigr].\leqno (2.8)$$

Fix $p>d+1$. Choose a constant  $T_o>0$ small enough such that $2pe^{-C_pT_0}>d+1$. It follows from (2.8) and Kolmogorov's  modification theorem that there exists  a version of $X_t(w,x_0)$, denoted by $\tilde{X}_t(w,x_0)$, such that $(t,x_o)\ra \tilde{X}_t(w,x_o)$ is continuous over $[0, T_o]\times \{|x_o|\leq R+1\}$ almost surely. But
$$\dis X_t(x_o,w)=x_o\quad\hbox{\rm if}\quad |x_o|>R.$$

We conclude that $(t,x_o)\ra \tilde{X}_t(x_o,w)$ can be extended continuously to
$[0, T_o]\times \R^d$. Let $\dis (\theta_{T_o}w)(t)=w(t+T_o)-w(T_o)$. Define for $0<t\leq T_o$,

$$\dis \tilde{X}_{T_o+t}(x_o,w)=\tilde{X}_t\bigl(\tilde X_{T_o}(x_o,w),\theta_{T_o}w\bigr).$$
Then $\dis  \tilde{X}_{T_o+\cdot}(x_o,w)$ satisfies the s.d.e $(0.1)$ driven by the Brownian motion $\theta_{T_o}w$ with the initial condition $\tilde X_{T_o}(x_o,w)$. By pathwise uniqueness, 
we see that  $\dis  \tilde{X}_{T_o+t}(x_o,w)={X}_{T_o+t}(x_o,w)$ almost surely for all 
$t\in [0, T_o]$. This means that $\tilde X_t(x_o,w)$ is a continuous version of $X_t(x_o,w)$ over $|0,2T_o]\times \R^d$. Continuing in this way, we get a continuous version 
on the whole space $[0, +\infty[\times \R^d$.

\medskip
{\bf Step 2:} general case.
\medskip
For $R>0$, let $f_R(x)$ denote a smooth function with compact support satisfying
$$\dis f_R(x)=1\quad\hbox{\rm for}\ |x|\leq R \quad\hbox{\rm and}\quad f_R(x)=0\quad\hbox{\rm for}\ |x|> R+1.
 \leqno(2.9)$$
Define 
$$\dis \sigma_R(x)=\sigma (x)f_R(x) \quad\hbox{\rm and}\quad f_R(x)=b(x)f_R(x). 
\leqno(2.10)$$
Let $X_t^R(x,w)$ be the unique solution of the s.d.e. (0.1) with $\sigma$ and $b$ replaced by 
$\sigma_R$ and $b_R$. Let $\tilde{X}_t^R(x,w)$ denote a continuous version of $X_t^R(x,w)$. Such a version exists according to step 1.    For $K>0$, set

$$\dis \tau_K^R(x)=\inf \{t>0;\ |\tilde{X}_t^R(x,w)|\geq K\}
\leqno(2.11)$$
$$\dis \tau_K(x)=\inf \{t>0;\ |X_t(x,w)|\geq K\}
\leqno(2.12)$$
By the pathwise uniqueness, for $|x|\leq K$, we have 
$$\dis \tau_K(x)=\tau_K^{K+2}(x)
\leqno(2.13)$$
almost surely.
For $|x|\leq R$, define
$$\dis \tilde{X}_.(x,w)=\tilde{X}_.^{R+2}(x,w)\quad \hbox{\rm on} \quad [0,\tau_R^{R+2}(x))\leqno (2.14)$$
Then it is clear that $\tilde{X}_.(x,w)$ is a version of ${X}_.(x,w)$. Let us prove that $\tilde{X}_t(x,w)$ is continuous in $(t,x)$ for almost all  $w$. Fix $x_0$ with $|x_0|\leq K$. Since the life time of the solution is infinity, there exists $R>0$ such that 
$\tau_R^{R+2}(x_0)>t$. This implies that
$\sup_{0\leq s\leq t}|\tilde{X}_s^{R+2}(x_o,w)|<R$. By the continuity, we can find a  neighbourhood $B_{\delta}(x_o)$ of $x_o$ such that $\sup_{0\leq s\leq t}|\tilde{X}_s^{R+2}(x,w)|<R$ or $ \tau_R^{R+2}(x)>t$ for all $x\in B_{\delta}(x_o)$. Hence, $\tilde{X}_s(x,w)=\tilde{X}_s^{R+2}(x,w)$ for all $x\in B_{\delta}(x_o)$ and $s\leq t$, which implies that $\tilde{X}_s(x_o,w)$ is continuous with respect to $(s,x_o)$.
\fin

\medskip
{\bf Remark 2.4:} Consider $\dis r(s)=\log{1\over s}\cdot\log{\log{1\over s}}$ for $s\in ]0, 1/2e]$. Clearly $s\ra r(s)$ is decreasing and $s\ra sr(s)$ is concave over $]0, 1/2e]$. Applying $(2.7)$, we get
$$\dis \varphi(t) \leq \exp{\Bigl(-\bigl[\log{1\over\varphi(0)}\bigr]^{e^{-C_p t}}\Bigr)}.$$
In order to apply the Kolmogorov's  modification theorem, we have to find $\alpha>0$ such that
$$\dis \exp{\Bigl(-\bigl[\log{1\over\varphi(0)}\bigr]^{e^{-C_p t}}\Bigr)}
\leq \varphi(0)^\alpha,$$
or
$$\dis \bigl[\log{1\over\varphi(0)}\bigr]^{e^{-C_p t}}
\geq \alpha\log{1\over\varphi(0)}$$
which is impossible when $|x_o-y_o|$ is small for any $t>0$.\fin

\vskip 8mm
{\bf 3. Flow of homeomorphisms}
\bigskip
{\bf Theorem 3.1}\quad {\it Assume that for $|x-y|\leq {1\over 2}$

$$\dis
\left\{\matrix{||\sigma(x)-\sigma(y)||^2&\leq& C\, |x-y|^2\,\log{1\over |x-y|},\cr
|b(x)-b(y)|&\leq& C\,  |x-y|\,\log{1\over |x-y|},}\right.$$
and $\sigma$ and $b$ are compactly supported. Then the solution of the s.d.e $(0.1)$ admits 
a version $X_t(x_o,w)$ such that 
$x_o\ra X_t(x_o,w)$ is a homeomorphism of $\R^d$ almost surely for all $t>0$.}
\medskip

{\bf Proof.} Let $R>0$ such that
$$\dis \sigma(x)=0\quad\hbox{\rm and}\quad
b(x)=0\quad\hbox{\rm for }\ |x|\geq R.$$

By section 2, the solution of s.d.e. $(0.1)$ admits a continuous version, still denoted by
$X_t(x_o,w)$, in $(t,x_o)$, such that $|X_t(x_o,w)|\leq R+1$ if $|x_o|\leq R+1$. So 
$x_o\ra \tilde X_t(x_o,w)$ defines a continuous map from $B(R+1)$ to $B(R+1)$, where
$$\dis B(r)=\bigl\{x\in\R^d;\ |x|\leq r\bigr\}.$$

\medskip
{\bf Lemma 3.2}\quad {\it Let $x_o\neq y_o$ and $\alpha<0$. Then there exists $C_\alpha, K_\alpha>0$ such that
$$\dis \E\Bigl(|X_t(x_o)-X_t(y_o)|^{2\alpha}\Bigr)
\leq C_\alpha\, |x_o-y_o|^{2\alpha e^{-K_\alpha\, t}},
\quad\hbox{\rm for}\ x_o, y_o\in B(R+1).
\leqno(3.1)$$}

\medskip
{\bf Proof.}\quad Let $\dis 0<\ge <{|x_o-y_o|^2\over 8(R+1)^2}$. Consider
$\dis \eta_t(w)=X_t(x_o)-X_t(y_o)$ and
$$\dis \xi_t(w)={|\eta_t(w)|^2\over 8(R+1)^2}.$$
 By theorem 1.1, we know that $\xi_t(w)\neq 0$ almost surely for all $t>0$. Define

$$\dis \tau_\ge =\inf\bigl\{t>0,\ \xi\leq \ge\bigr\}.$$

By It\^o formula, for $s<t$,
$$\dis \xi_{t\wedge\tau_\ge}^\alpha-\xi_{s\wedge\tau_\ge}^\alpha
=\int_{s\wedge\tau_\ge}^{t\wedge\tau_\ge}\alpha \xi_u^{\alpha-1}\, d\xi_u
+{1\over 2}\alpha(\alpha-1)\int_{s\wedge\tau_\ge}^{t\wedge\tau_\ge}\alpha \xi_u^{\alpha-2}\, d\xi_u\cdot d\xi_u.$$

Where
$$\dis\eqalign{ d\xi_t ={1\over 4(R+1)^2}\Bigl\{ 
&\bigl<\eta_t(w), \bigl(\sigma(X_t(x_o))-\sigma(X_t(y_o))\bigr) dw_t\bigr>\cr
&+\bigl<\eta_t(w), b(X_t(x_o))-b(X_t(y_o))\bigr>\, dt\cr
&+ {1\over 2}||\sigma(X_t(x_o))-\sigma(X_t(y_o))||^2\, dt\Bigr\},}$$
and
$$\dis d\xi_t\cdot d\xi_t = \Bigl({1\over 4(R+1)^2}\Bigr)^2
|\bigl(\sigma^*(X_t(x_o))-\sigma^*(X_t(y_o))\bigr)\eta_t|^2\, dt.$$

Using lemma 2.1 for $p=1$, there exists a constant $C>0$ such that
$$\dis |b(X_t(x_o))-b(X_t(y_o))|\leq C\, |\eta_t|\log{8(R+1)^2\over |\eta_t|^2 }
\leqno(3.2)$$
 and
$$\dis ||\sigma(X_t(x_o))-\sigma(X_t(y_o))||^2
\leq C\, |\eta_t|^2\log{8(R+1)^2\over |\eta_t|^2 }.
\leqno(3.3)$$

Therefore by $(3.2)$,
$$\dis\eqalign{& \Bigl|{\alpha \xi_u^{\alpha-1}}\bigl<\eta_u, b(X_t(x_o))-b(X_t(y_o))\bigr>\Bigr|\cr
&\leq -2C\alpha\, \xi_u^\alpha\log{1\over\xi_u}
=-2C\, \xi_u^\alpha\log{1\over\xi_u^\alpha}.}$$

Similarly , we have 
$$\dis\eqalign{
&{1\over 2}\alpha(\alpha-1)\xi_u^{\alpha-2}\Bigl({1\over 4(R+1)^2}\Bigr)^2
|\bigl(\sigma^*(X_t(x_o))-\sigma^*(X_t(y_o))\bigr)\eta_t|^2\cr
&\hskip 10mm\leq 2C|\alpha-1|\,\xi_u^\alpha\log{1\over\xi_u^\alpha}.}$$

Let $dM_t$ be the martingale part of $d\xi_t$. Then we have
$$\dis \xi_{t\wedge\tau_\ge}^\alpha-\xi_{s\wedge\tau_\ge}^\alpha
\leq\int_{s\wedge\tau_\ge}^{t\wedge\tau_\ge}\alpha \xi_u^{\alpha-1}\, dM_u
+K(\alpha)\int_{s\wedge\tau_\ge}^{t\wedge\tau_\ge} \xi_u^{\alpha}\log{1\over \xi_u^\alpha}\,du .$$
Define $\dis \varphi(t)=\E(\xi_{t\wedge\tau_\ge}^\alpha)$. From the above inequality, we get
$$\dis \varphi'(t)\leq K(\alpha)\varphi(t)\log{1\over\varphi(t)},$$
which implies that
$$\dis \varphi(t)\leq \varphi(0)^{e^{-K(\alpha)\, t}}.$$

Hence,
$$\dis \E\bigl(\xi_{t\wedge\tau_\ge}^\alpha)
\leq \Bigl({|x_o-y_o|^2\over 8(R+1)^2}\Bigr)^{\alpha\,e^{-K(\alpha)\, t}}
\quad\hbox{\rm for}\ x_o,y_o\in B(R+1).$$
Since $\tau_\ge\uparrow +\infty$ while $\ge\downarrow 0$, we get
$$\dis \E\bigl(\xi_{t}^\alpha)
\leq \Bigl({|x_o-y_o|^2\over 8(R+1)^2}\Bigr)^{\alpha\,e^{-K(\alpha)\, t}}
\quad\hbox{\rm for}\ x_o,y_o\in B(R+1)$$
or $(3.1)$ holds.\fin

\medskip
{\bf Lemma 3.3}\quad {\it Let $\delta>0$ and 
$\dis\Delta_\delta
=\bigl\{(x_o,y_o)\in B(R+1)\times B(R+1);\ |x_o-y_o|\geq\delta\bigr\}$.
Take a continuous version $X_t(x_o)$ and set
 $\dis \eta_t(x_o,y_o)=|X_t(x_o)-X_t(y_o)|^{-1}$. Then  for any $p>1$, there exist constant $\dis C_{p,T,\delta}, K_1(p)>0$ such that
$$\dis \eqalign{&\E\Bigl( |\eta_t(x_o,y_o)-\eta_s(\tilde{x}_o,\tilde{y}_o)|^p\Bigr)\cr
&\leq C_{p,T,\delta}\Bigl( |x_o-\tilde x_o|^{pe^{-K_1(p)T}}
+|y_o-\tilde y_o|^{pe^{-K_1(p)T}}+|t-s|^{p/2}\Bigr)}
\leqno(3.4)$$
for
$\dis (x_o,y_o), (\tilde x_o,\tilde y_o)\in \Delta_\delta$ and $s,t\in [0,T]$.}

\medskip
{\bf Proof.} We have
$$\dis\eqalign{& |\eta_t(x_o,y_o)-\eta_s(\tilde x_o,\tilde y_o)|\cr
&\hskip -10mm
=\eta_t(x_o,y_o)\eta_s(\tilde x_o,\tilde y_o)\Bigl| |X_t(x_o)-X_t(y_o)|
-|X_s(\tilde x_o)-X_s(\tilde y_o)|\Bigr|\cr
&\hskip -10mm
\leq \eta_t(x_o,y_o)\eta_s(\tilde x_o,\tilde y_o)\Bigl( |X_t(x_o)-X_s(\tilde x_o)|
+|X_t( y_o)-X_s(\tilde y_o)|\Bigr).}$$

By lemma 3.1, there exists a constant $C_{p,T,\delta}>0$ such that
$$\dis \E\bigl(\eta_t(x_o,y_o)^{4p}\bigr)\leq C_{p,T,\delta}
\quad\hbox{\rm for}\ (x_o,y_o)\in \Delta_\delta,\, t\in[0,T].$$
Now combining with $(2.8)$, we get 
$$\dis \eqalign{&\E\Bigl( |\eta_t(x_o,y_o)-\eta_s(\tilde{x}_o,\tilde{y}_o)|^p\Bigr)\cr
&\leq C_{p,T,\delta}\Bigl( |x_o-\tilde x_o|^{2pe^{-K_1(p)T}}
+|y_o-\tilde y_o|^{2pe^{-K_1(p)T}}+2|t-s|^{p}\Bigr)^{1/2}}$$
which is dominated by the right hand side of $(3.4)$. So we get the result.\fin

\medskip
{\bf Proof of theorem 3.1}\quad By $(3.4)$, for $p>2(2d+1)$ and $T_o>0$ small enough, we can apply the Kolmogorov's modification theorem to get that  $\eta_t(x_o,y_o)$ has a continuous version 
$\dis\tilde \eta_t(x_o,y_o)$ on $[0,T_o]\times\Delta_o$. This means that 
$\dis (t,x_o,y_o)\ra \tilde \eta_t(x_o,y_o)$ is continuous on $[0,T_o]\times\Delta_o$ almost surely. Let
$D$ be a countable dense subset of $[0,T_o]\times\Delta_o$. Then almost surely, 
for all $(t,x_o,y_o)\in D$,
$$\dis \tilde\eta_t(x_o,y_o)=|X_t(x_o)-X_t(y_o)|^{-1},$$
or 
$$\dis |X_t(x_o)-X_t(y_o)|=\tilde\eta_t(x_o,y_o)^{-1}.
\leqno(3.5)$$

Now by continuity, the relation $(3.5)$ holds for all $(t,x_o,y_o)\in [0,T_o]\times\Delta_o$.
In particular, this relation shows that almost surely for all $t>0$, $x_o\ra X_t(x_o)$ is injective on $B(R+1)$. Since $X_t(x_o)=x_o$ for $|x_o|>R$, $x_o\ra X_t(x_o)$ is injective on
the whole $\R^d$.  Let $\dis\overline{\R^d}=\R^d\cup\{\infty\}$ which is homeomorph to the sphere $S^d$. Extend the map $X_t$ to $\dis\overline{\R^d}$ by setting
$\dis X_t(\infty)=\infty$. It is clear that $X_t$ is continuous near $\infty$. So $X_t$ can be seen as a continuous map from $S^d$ into $S^d$. Since $X_t$ is homotope to the identity, we conclude that $X_t$ is a surjective map on $S^d$. Therefore $X_t$ is a homeomorphism of $S^d$ and  its restriction $X_t$ on $\R^d$ is a homeomorphism of $\R^d$. Now using the relation
$$\dis X_{T_o+t}(x_o,w)=X_t\bigl(X_{T_o}(x_o,w), \theta_{T_o}w\bigr),$$
we conclude that $X_t$ is a homeomorphsim of $\R^d$ for all $t>0$.\fin
\vskip 0.3cm
{\bf Acknowledgements:} The owrk of T.S.Zhang is partially supported by the British EPSRC (grant no. GR/R91144/01)

\vskip 10mm
\centerline{\bf References}

\bigskip
[El] Elworthy K.D., Stochastic flows on Riemannian manifolds, in {\it Diffusion processes and related problems in analysis}, ed. by M.A. Pinsky and V. Wihstutz, Birkh\"auser Boston, 1992.

[Em] Emery M.: Non confluence des solutions d'une equation stochastique lipschitzienne. Seminaire Proba. XV. Lecture Notes in Mathematics, vol. 850 (587-589) Pringer, Berlin Heidelberg New York 1981.

[Fa] Fang S.: Canonical Brownian motion on the diffeomorphism group of the circle, J. Funct. Anal. {\bf 196} (2002), 162-179.

[FZ1] Fang S. and Zhang T.S. : Stochastic differential equations with non-Lipschitz coefficients: pathwise uniqueness and no explosion. To appear in CRAS...

[FZ2] Fang S. and Zhang T.S. : Study of Stochastic differential equations with non-Lipschitz coefficients: I. pathwise uniqueness and large deviations. Preprint 2003.

[IW] Ikeda I. , Watanabe S.: ``Stochastic differential equations and Diffusion processes,'' North-Holland, Amsterdam, 1981.

[Ku] Kunita H. : ``Stochastic flows and stochastic differential equations'' Cambridge University Press 1990.

[LJR1] Le Jan Y. , Raimond O.: Integration of Brownian vector fields, Annals of Prob.
 {\bf 30}  (2002),  no. 2, 826-873.

[LJR2] Le Jan Y. , Raimond O.: Flows, coalescence and noise, Annals of Prob. 2003.

[Ma1] Malliavin P.: The Canonical diffusion above the diffeomorphism group of the circle, C.R. Acad. Sci. Paris, S\'erie I {\bf 329} (1999), 325-329.

[Ma2] Malliavin P.: {\it Stochastic Analysis}, Grunlehren des Math. {\bf 313}, Springer, 1997.

[Pr] Protter P. : ``Stochastic integration and differential equations'' Springer-Verlag Berlin Heidelberg New York 1990.

[RY] Revuz D., Yor M.: ``Continuous martingales and Brownian motion'', Grund. der Math. Wissenschaften {\bf 293}, 1991, Springer-Verlag.

[YO] Yamada T. and Ogura Y. : On the strong comparison theorems for solutions of stochastic differential equations, Z. Wahrscheinlichkeitstheorie verw. Gebiete 56 (1981) 3-19.
\end